\documentclass[12pt]{article}%
\usepackage{amsmath}
\usepackage{amsfonts}
\usepackage{amssymb}
\usepackage{graphicx}%
\providecommand{\U}[1]{\protect\rule{.1in}{.1in}}
\newtheorem{theorem}{Theorem}

\newtheorem{example}[theorem]{Example}

\newtheorem{proposition}[theorem]{Proposition}
\newtheorem{remark}[theorem]{Remark}

\newenvironment{proof}[1][Proof]{\noindent\textbf{#1.} }{\ \rule{0.5em}{0.5em}}
\begin{document}

\title{A counterexample regarding the relatively uniform completion of a principal ideal.}
\author{Youssef Azouzi\thanks{The authors are members of the GOSAEF research group}
\and {\small Research Laboratory of Algebra, Topology, Arithmetic, and Order}\\{\small Department of Mathematics}\\COSAEF {\small Faculty of Mathematical, Physical and Natural Sciences of
Tunis}\\{\small Tunis-El Manar University, 2092-El Manar, Tunisia}}
\maketitle

\begin{abstract}
We present a counterexample related to relative uniform convergence, showing
that, in general, the relatve uniform completion of the principal ideal of a
vector lattice $E$ generated by an element $x$ is stricly contained in the
ideal generated by $x$ in the relatively uniform completion of $E.$

\end{abstract}

\section{Introduction}

In the theory of vector and Banach lattices, various types of convergence are
considered, including both topological and non-topological types. Recently,
the systematic study of convergence structures has gained attention among
researchers (see, for example, \cite{L-748} and \cite{L-999}). Among these is
the concept of relatively uniform convergence. However, notable observation
emerges from the literature: there lacks a systematic and comprehensive study
of this type of convergence. Moreover, fundamental and basic questions
regarding relatively uniform convergence remain unanswered, posing significant
challenges. In this short note we will present an example addressing one of
these elusive questions. It compares the relatively uniform completion of the
principal ideal of a vector lattice $E,$ generated by $x$ with the principal
ideal generated by the same element in the relatively uniform completion of
the whole space. Let us first introduce some notation and terminology, and
review a few basic notions related to our subject. Let $E$ be an Archimedean
vector lattice. A sequence $\left(  x_{n}\right)  $ in $E$ is said to be
relatively uniformly convergent to $x$ if there exists (a regulator) $u\in
E^{+}$ and a sequence $\left(  \varepsilon_{n}\right)  $ of reals convergent
to $0$ such that $\left\vert x_{n}-x\right\vert \leq\varepsilon_{n}u$ for all
$n\in\mathbb{N}.$ If $E$ is Archimedean then relatively uniform limits are
unique. We say that $\left(  x_{n}\right)  $ is relatively uniform Cauchy if
there exists a positive element $u$ in $E$ and a real sequence $\left(
\varepsilon_{n}\right)  $ convergent to $0$ such that $\left\vert x_{n}%
-x_{m}\right\vert \leq\varepsilon_{n}u$ for all $n,m\in\mathbb{N}$ with $m\geq
n.$ The vector lattice $E$ is relatively uniformly complete if every
relatively uniform Cauchy sequence in $E$ is relatively uniformly convergent.
Another way to describe uniform completeness is the folllowing: For each
element $e$ in $E^{+},$ the ideal $E_{e}$ generated by $e$ can be equipped
with the following norm%
\[
\left\Vert x\right\Vert _{e}=\inf\left\{  \lambda\geq0:\left\vert x\right\vert
\leq\lambda e\right\}  .
\]

Then the space $\left(  E_{e},\left\Vert .\right\Vert _{e}\right)  $ is a
normed lattice and if $E$ is a Banach lattice with strong unit $e,$ then
$\left\Vert .\right\Vert _{e}$ is equivalent to the original norm. It is not
difficult to show that $E$ is uniformly complete if and only if $\left(
E_{e},\left\Vert .\right\Vert _{e}\right)  $ is a Banach lattice for every
$e\in E^{+}.$ Banach lattices and Dedekind $\sigma$-complete vector lattices
are uniformly complete. We denote by $E^{\delta}$ the Dedekind completion of
$E$ and by $E^{ru}$ its relatively uniform completion. It is a well-known fact
that $E^{ru}$ is the intersection of all uniformly complete vector sublattices
of $E^{\delta}$ that contain $E.$ The ideal generated by $x\in E$ will be
denoted by $E_{x}$. Every Dedekind complete vector lattice is relatively
uniformly complete. It is easy to see that if $x\in E$ then $\left(
E^{\delta}\right)  _{x}$ can be naturally identified with $\left(
E_{x}\right)  ^{\delta}.$ However, when considering the relative uniform
completion, the situation becomes more subtle and complex. In this case, only
the inclusion $\left(  E_{x}\right)  ^{ru}\subseteq\left(  E^{ru}\right)
_{x}$ is obvious. We provide in this note an example showing that this
inclusionmay be strict. For more details about relatively uniform convergence
the reader is referred to \cite{b-241}.

\section{A counter-example}

It is clear that if $E$ is an Archimedean vector lattice and $v$ is an element
of $E,$ then the uniform completion of the ideal generated by $v$ is contained
in the ideal generated by $v$ in the uniform completion of $E,$ that is,%
\[
\left(  E_{v}\right)  ^{ru}\subseteq\left(  E^{ru}\right)  _{v}.
\]
In a recent work, the authors of \cite{L-1045} claimed that this inclusion is
actually an equality. However, the following example shows that the inclusion
can, in fact, be strict. This mistaken claim is central to their proof of a
surprising result \cite[Theorem 2]{L-1045}, namely that every element in the
r.u-completion of an Archimedean vector lattice $E$ is the r.u-limit of a
sequence in $E.$ This statement is, however, known to be false, even though no
explicit counterexample appears to be available in the literature. Although
Quinn in \cite{L-1005} presented several counterexamples, he claimed that he
was not able to find one specifically adressing this particumlar isuue.
Nevertheless Ball and Hager have stated that such a counterexample does exist
\cite{BH}.

\begin{example}
\label{E1}Consider the vector lattice $E$ consisting of all linear piecewise
continuous functions on $\left[  0,1\right]  .$ Let $u\in E$ be defined by
$u\left(  t\right)  =t$ for all $t\in\left[  0,1\right]  .$ We have
$E^{ru}=C\left[  0,1\right]  $ and then $\left(  E^{ru}\right)  _{u}$ is the
principal ideal of $C\left[  0,1\right]  $ generated by $u.$ We claim that%
\[
\left(  E_{u}\right)  ^{ru}=\left\{  f\in E:\exists\varphi_{n}\in
E_{u},\varepsilon_{n}>0:\left\vert \varphi_{n}-f\right\vert \leq
\varepsilon_{n}u,\varepsilon_{n}\longrightarrow0\right\}  .
\]
Furthermore the inclusion $\left(  E_{u}\right)  ^{ru}\subset\left(
E^{ru}\right)  _{u}$ is strict.
\end{example}

\begin{proof}
Denote by $H$ the left side of the above equality. Then the inclusion
$H\subseteq\left(  E_{u}\right)  ^{ru}$ is obvious. It is enough so to prove
that $H$ is relatively uniformly complete. Assume that $\left(  f_{n}\right)
$ is an r.u-Cauchy sequence in $H.$ Then there exist a real sequence $\left(
\varepsilon_{n}\right)  $ converging to $0$ and an element $v$ in $\left(
E_{u}\right)  ^{ru}$ such that%
\begin{equation}
\left\vert f_{n}-f_{m}\right\vert \leq\varepsilon_{n}v,\qquad\text{for }m\geq
n\geq1.\label{T1}%
\end{equation}
Now as $\left(  E_{u}\right)  ^{ru}$ is contained in $\left(  E^{ru}\right)
_{u}$ we may assume that $v=u.$ Let $f_{n}$ denote the uniform limit of
$\left(  f_{n}\right)  $ in $E.$ By letting $m\longrightarrow\infty$ in
(\ref{T1}) we get%
\begin{equation}
\left\vert f_{n}-f\right\vert \leq\varepsilon_{n}u,\qquad\text{for }%
n\geq1.\label{T2}%
\end{equation}
Moreover as $f_{n}\in H$ for each $n$ we can find $\varphi_{n}\in E_{u}$ such
$\left\vert \varphi_{n}-f_{n}\right\vert \leq\dfrac{1}{n}u.$ Hence the
following inequality holds%
\[
\left\vert \varphi_{n}-f\right\vert \leq\left(  \dfrac{1}{n}+\varepsilon
_{n}\right)  u,
\]
for all $n\in\mathbb{N}.$ This shows that $f\in H.$ Also by (\ref{T2}),
$f_{n}\overset{r.u.}{\longrightarrow}f$ in $H.$ This shows that $H$ is
relatively uniformly complete as it was claimed.

Next we will show that $\left(  E_{u}\right)  ^{ru}$ is contained strictly in
$\left(  E^{ru}\right)  _{u}.$ To this end consider tow real sequences
$\left(  a_{n}\right)  $ and $\left(  b_{n}\right)  $ satisfying the following
two conditions:

(i) $0<a_{n+1}<b_{n+1}<a_{n}<b_{n}<....<a_{1}<b_{1}=1,$ and

(ii) $\lim b_{n}=0=\lim a_{n}.$

Then define a function $f$ on $\left[  0,1\right]  $ by putting $f\left(
0\right)  =0,$ $f\left(  b_{n}\right)  =0,$ $f\left(  a_{n}\right)  =a_{n}$
for $n\in\mathbb{N},$ and $f$ is linear on each of the intervals $\left[
a_{n},b_{n}\right]  ,\left[  b_{n},a_{n-1}\right]  $. It is easy to check that
$f$ is continuous and $0\leq f\leq u.$ Hence $f\in$ $\left(  E^{ru}\right)
_{u}$. We claim that $f$ does not belong to$\left(  E_{u}\right)  ^{ru}=H.$ We
will argue by contradiction in assuming that $f\in\left(  E_{u}\right)  ^{ru}$
and then it satisfies, by the first part, the following conditions%
\[
\left\vert \varphi_{n}-f\right\vert \leq\varepsilon_{n}u\text{ for all }%
n\in\mathbb{N}.
\]
for some sequences $\left(  \varphi_{n}\right)  $ in $E$ and $\left(
\varepsilon_{n}\right)  $ in $\left(  0,\infty\right)  $ with $\varepsilon
_{n}\downarrow0.$ Pick an integer $n_{0}$ such that $\varepsilon_{n_{0}}%
\leq1/3.$ As $\varphi_{n_{0}}\in E_{u}$ there exist $\delta>0,\lambda
\in\mathbb{R}$ such that $\varphi_{n_{0}}\left(  t\right)  =\lambda t$ for
$t\in\left[  0,\delta\right]  .$ Thus%
\[
\left\vert \lambda t-f\left(  t\right)  \right\vert \leq1/3.t\text{ for }%
t\in\left[  0,\delta\right]  .
\]
Applying this to $t=b_{k}$ for $k$ large enough yields%
\[
\left\vert \lambda\right\vert \leq1/3.
\]
On the other hand if we take $t=a_{k}$ for $k$ large enough we get%
\[
\left\vert \lambda-1\right\vert \leq1/3.
\]
These last two inequalities are incompatible. This contradiction shows that
$f\notin\left(  E_{u}\right)  ^{ru}$ and then the proof is complete.
\end{proof}

\begin{remark}
Our example shows that the proof of Lemma 1 given in \cite{L-1045} is
incorrect. But sometimes it is not easy to identify the error. The authors
used the following equality%
\[
\bigcap\limits_{\substack{F\text{ closed in }E^{\delta}\\E\subseteq F}%
}F\cap\left(  E^{\delta}\right)  _{x}=\bigcap\limits_{\substack{G\subseteq
E^{\delta}\text{ closed}\\E_{x}\subseteq G}}G
\]
which does not hold in general. While the left-hand side is indeed contained
in the right-hand side, the converse inclusion fails. Not every closed
sublattice $G$ in $E^{\delta}$ containing $E_{x}$ can be written in the form
$G=F\cap\left(  E^{\delta}\right)  _{x}$ for some closed vector sublattice $F$
with $E\subset F\subset E^{\delta}.$
\end{remark}

\begin{remark}
It was asserted in \cite{L-1045} that if $E$ \textit{is a vector lattice, then
for every} $x\in E^{ru}$, \textit{there exists a sequence} $(x_{n})$ in $E$
\textit{such that} $x_{n}\overset{r.u.}{\longrightarrow}x.$ However, as was
noted above, his claim contradicts a well-known fact in the literature (see
for example \cite{b-241,L-1044}). Regarding the equality beteween the
r.u-closure and the r.u-completion of a vector lattice we propose the
following "positive" result. While it aligns with the intended direction of
the authors of \cite{L-1045}, it does not directly contribute to achieving
their final aim. Denote for a vector lattice the ru-colsure of $E$ in
$E^{\delta}$ by $\overline{E}^{ru}.$
\end{remark}

\begin{proposition}
Let $E$ be a Riesz space such that for each $x\in E$ we have $E_{x}%
^{ru}=\left(  E^{ru}\right)  _{x}.$ Then $E=\overline{E}^{ru}.$ That is every
element in $E^{ru}$ is an ru limiot of a sequence in $E.$
\end{proposition}

\end{document}